\documentclass[a4paper,11pt]{article}

\usepackage{url}
\usepackage[latin1]{inputenc}
\usepackage[english]{babel}
\usepackage{graphicx}
\usepackage{amssymb}
\usepackage{amsmath}
\usepackage{amsfonts}
\usepackage{algorithm}
\usepackage{algorithmic}
\usepackage{amsthm}

\usepackage{setspace}

\newtheorem{theorem}{Theorem}
\newtheorem{lemma}{Lemma}
\newtheorem{corollary}{Corollary}
\theoremstyle{definition} 
\newtheorem{example}{Example}
\newtheorem{remark}{Remark}
\DeclareMathOperator{\argmin}{argmin}
\def \matlab    {MATLAB$^{\textrm{\tiny \textregistered}}$}

\begin{document}
\renewcommand{\labelitemi}{$\diamond$}


\title{Using Underapproximations \\ for Sparse Nonnegative Matrix Factorization}

\author{Nicolas Gillis \and Fran\c{c}ois Glineur\thanks{Center for Operations Research and Econometrics, Universit\'e catholique de Louvain, Voie du Roman Pays, 34, B-1348 Louvain-La-Neuve, Belgium ;
nicolas.gillis@uclouvain.be and francois.\mbox{glineur}@uclouvain.be. The first author is a research fellow of the Fonds de la Recherche Scientifique (F.R.S.-FNRS). This text presents research results of the Belgian Program on Interuniversity Poles of Attraction initiated by the Belgian State, Prime Minister's Office, Science Policy Programming. The scientific responsibility is assumed by the authors.} }

\date{\small October 2009}


\maketitle
\thispagestyle{empty}

\begin{abstract}
Nonnegative Matrix Factorization consists in (approximately) factorizing a nonnegative data matrix by the product of two low-rank nonnegative matrices. It has been successfully applied as a data analysis technique in numerous domains, e.g., text mining, image processing, microarray data analysis, collaborative filtering, etc.

We introduce a novel approach to solve NMF problems, based on the use of an underapproximation technique, and show its effectiveness to obtain sparse solutions. This approach, based on Lagrangian relaxation, allows the resolution of NMF problems in a recursive fashion. We also prove that the underapproximation problem is NP-hard for any fixed factorization rank, using a reduction of the maximum edge biclique problem in bipartite graphs.

We test two variants of our underapproximation approach on several standard image datasets and show that they provide sparse part-based representations with low reconstruction error. Our results are comparable and sometimes superior to those obtained by two standard Sparse Nonnegative Matrix Factorization techniques.\\
\bigskip

\end{abstract}

\newpage

\section{Introduction} \label{intro}

Nonnegative Matrix Factorization (NMF) is a recent data analysis technique with applications in image processing, text mining, spectral unmixing, air emission control, computational biology, clustering, etc. (see \cite{Ber, Dhi, Dev, diep} and references therein). NMF can be described as follows: given a nonnegative input matrix $M \in \mathbb{R}^{m \times n}_+$ and an integer $1 \le r < \min(m,n)$, find two nonnegative matrices $V \in \mathbb{R}^{m \times r}_+$ and $W \in \mathbb{R}^{r \times n}_+$ whose product approximates the input matrix as closely as possible:
\begin{equation}
M \approx V W ,
\label{approx}
\end{equation}
so that $VW$ is a low-rank approximation of $M$. Matrix factorization can be interpreted as a linear factor model: assuming that each column of the input matrix $M$ represents an element of a data set, decomposition~\eqref{approx} can be written as\footnote{In this text, $A_{ij}$ stands for the $(i,j)$-entry of a matrix $A$, $A_{:j}$ for its $j^{\textrm{th}}$ column and $A_{i:}$ for its $i^{\textrm{th}}$ row.}
\begin{equation}
M_{:j} \approx \sum_{k} V_{:k} W_{kj}, \quad \forall j\;, \nonumber
\end{equation}
i.e., each input column $M_{:j}$ is a linear combination of a set of $r$ basis elements $V_{:k}$ with corresponding weights $W_{kj}$.

\begin{figure*}[ht!]
\begin{center}
\includegraphics[width=\textwidth]{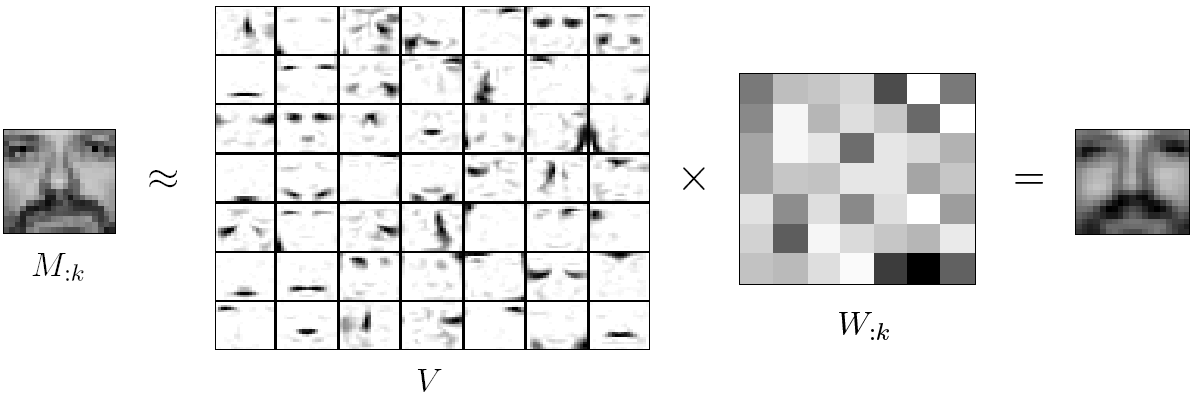}
\caption{NMF applied to the CBCL Face Database $\#1$, MIT Center For Biological and Computation Learning (available at $http$:$//cbcl.mit.edu/cbcl/software$-$datasets/FaceData2.html$). It consists of the approximation of 2429 gray-level images of faces represented with $19 \times 19$ pixels (columns of $M$) using $r=49$ basis elements (columns of $V$).}
\label{facto}
\end{center}
\end{figure*}

In contrast with standard linear factor model techniques such as Principal Component Analysis, NMF considers nonnegativity of the input columns to be an important feature 
and consequently requires the basis elements to be also nonnegative, so that they can be interpreted in the same way (e.g., these columns can correspond to images described by nonnegative pixel intensities or to texts represented by vectors of nonnegative word counts). Furthermore, NMF imposes nonnegativity of the weights, leading to an essentially additive reconstruction of the input columns by the basis elements. This representation is then \emph{part-based}: basis elements $V_{:k}$ will represent common parts of the columns of $M_{:j}$. For example, if each column of $M$ represents a face using pixel intensities, the basis elements generated by NMF can be facial features, such as eyes, noses and lips, as shown in Figure~\ref{facto}.

This low-rank approximation technique with nonnegativity constraints was introduced in 1994 by Paatero and Tapper~\cite{Paa} and started to be extensively studied after the publication of an article
by Lee and Seung~\cite{LS1} in 1999. Since an exact representation of the input matrix cannot be obtained in general, the quality of the approximation is measured by some criterion, typically the sum of the squares of the errors on the entries, which leads to the following minimization problem\footnote{$||A||_F = (\sum_{i,j} A_{ij}^2)^\frac12$ denotes the Frobenius norm of matrix $A$.}:
\begin{equation}
\min_{V \in \mathbb{R}^{m \times r}, W \in \mathbb{R}^{r \times n}} ||M-VW||_F^2 \tag{NMF} \text{ such that }
 V \geq 0 \text{ and }W \geq 0. \label{NMF} \end{equation}

An important feature of NMF is that its nonnegativity constraints typically induce \emph{sparse factors}, i.e., factors with relatively many zero entries. Intuitively, decomposition into parts requires the basis elements to be sparse, cf.\@ Figure~\ref{facto}. More formally, the reason for this behavior is that stationary points $(V,W)$ of NMF will be typically located at the boundary of the feasible domain $\mathbb{R}_+^{m \times r} \times \mathbb{R}_+^{r \times n}$, hence will feature zero components. This can be explained with the first-order optimality conditions: because the set of stationary points of a problem of the type
\[
\min_{x \in \mathbb{R}^n} \;\; f(x) \text{ such that } x \ge 0
\]
is given by the following expression (where $\nabla f(x)$ is the gradient of $f$)
\[
S_f = \{ x \in \mathbb{R}^n \; | \; x \geq 0, \nabla f(x) \geq 0 \textrm{ and } \mathbf{x_i [\nabla f(x)]_i = 0} \, \forall i \},
\]
some components of the solution can be expected to be equal to zero.

Sparsity of the factors is an important consideration in practice: in addition to reducing memory requirements to store the basis elements and their weights, sparsity improves interpretation of the factors, especially when dealing with classification/clustering problems, \rm{e.g.}, in text mining \cite{SBPP06} and computational biology \cite{GC05, Kim}. By contrast, unconstrained low-rank approximations such as Principal Component Analysis (PCA) do not naturally generate sparse factors (for that reason, low-rank approximations techniques with additional sparsity constraints have been recently introduced; this is referred to as Sparse Principal Component Analysis, Sparse PCA or SPCA, see, e.g.,\@ \cite{AE07} and references therein).

Although solutions of NMF typically display some level of sparsity, some applications require even sparser solutions, leading to variants of NMF called Sparse Nonnegative Matrix Factorization.
They are in general developed in two different ways: some authors define a priori a desired sparsity level and adapt the main iteration of their method in order to guarantee that the factors satisfy that level of sparsity throughout the application of the algorithm, see, e.g., \cite{Hoy, HS06}. Alternatively, a penalty term can be added to the objective function to prevent the algorithm from considering dense solutions, see~\cite{LH01}. In particular, it is well-known that $l_1$-norm penalty terms induce sparser solutions (see, e.g., \cite{KP08, Kim, Cic4}). More details about these techniques are given at the beginning of Section~\ref{appl}.

Unfortunately the advantages of NMF (part-based representation and sparsity) over PCA come at a certain price. First, because of the additional nonnegativity constraints, the approximation error of the input data for a given  factorization rank $r$ will always be higher for NMF than in the unconstrained case. Second, optimization problem~\eqref{NMF} is more difficult to solve than its unconstrained counterpart: while PCA problems can be solved in polynomial time (e.g., using a Singular Value Decomposition technique \cite{Gol}), NMF problems belong to the class of NP-hard problems, as recently showed by Vavasis \cite{Vav}. However, it should also be pointed out that these drawbacks (higher error, NP-hardness) are also present for competing techniques emphasizing sparsity, such as SPCA.

Because of its NP-hardness, practical algorithms cannot be expected to find provably optimal global solutions for \eqref{NMF} in a reasonable amount of time and aim instead at finding locally optimal solutions. Most methods start from some initial guess factors $(V,W)$ and improve them iteratively using nonlinear optimization schemes such as projected gradient methods~\cite{Lin}, Newton-like methods~\cite{Dhi2, Zdu}, (block-)coordinate descent (also called alternating nonnegative least squares -- NNLS)~\cite{Chen, Cic, KP08}, multiplicative updates~\cite{LS2}, etc. (see also~\cite{Ber, Dhi, Cic2, Hod} and references therein).

In this paper, we introduce a novel approach to solve NMF problem based on the use of an underapproximation technique and show its effectiveness to obtain sparse solutions. Section~\ref{nmuS} introduces our underapproximation problem, motivated by a recursive technique to solve NMF, studies the sparsity of its solutions and proves that it is NP-hard for any fixed factorization rank. Nevertheless, Section~\ref{Lag} describes an algorithm to solve it approximately using a technique based on Lagrangian relaxation. Finally, in the last section, we test this approach on several standard image datasets,
and show both qualitatively and quantitatively that it provides part-based and sparse representations that are comparable and sometimes superior to those obtained with standard Sparse Nonnegative Matrix Factorization techniques.

\section{Nonnegative Matrix Underapproximation} \label{nmuS}
\subsection{A recursive approach}

Finding a rank-one nonnegative matrix factorization, i.e., solving \eqref{NMF} with $r=1$ is notably easier than for higher factorization ranks: while the general problem is NP-hard, computing a globally optimal rank-one approximation can be done in polynomial time. More specifically, the first rank-one factor of the Singular Value Decomposition (SVD) of the input matrix is an optimal solution: indeed, the Perron-Frobenius theorem implies that the dominant left and right singular vectors of a nonnegative matrix are nonnegative, while the Eckart-Young theorem states that the outer product of these dominant singular vectors is the best possible rank-one approximation in the Frobenius norm.

In principle, we might try exploit this result to find factorizations of higher ranks by applying it recursively: after identification of an optimal rank-one NMF solution $(v, w)$, one could subtract the $v w$ factor from $M$ and apply the same technique to $M - v w$ to recover the next rank-one factor. Unfortunately, this idea cannot work: the difference between $M$ and its rank-one approximation may contain negative values (typically roughly half of them), so that the next SVD factor will no longer provide a nonnegative solution. Moreover, there is no hope of replacing SVD by another efficient technique for this step since \cite{GG08} shows that it is NP-hard to find the optimal nonnegative rank-one approximation to a matrix which is not nonnegative.

If we wish to keep the principle of a recursive algorithm finding one rank-one factor at a time, we have to add a constraint ensuring that the $v w$ factor, when subtracted from $M$, gives a nonnegative remainder, i.e., we need to have $v w \le M$. Therefore we introduce a similar \emph{upper bound constraint} $VW \leq M$ to the general (NMF) problem and obtain a new problem we call \emph{Nonnegative Matrix Underapproximation} (NMU): given $M \in \mathbb{R}^{m \times n}_+$and $1 \leq r < \min(m,n)$, the NMU optimization problem is defined as
\begin{equation} \min_{V \in \mathbb{R}^{m \times r}, W \in \mathbb{R}^{r \times n}} ||M-VW||_F^2 \text{ such that } V \geq 0,\ W \geq 0 \text{ and } VW \leq M. \tag{NMU} \label{NMU}
\end{equation}
Assuming we are able to solve it for $r=1$, an underapproximation of any rank can then be built by following
the recursive procedure outlined above. More precisely, if $(V_{:1},W_{1:})$ is a rank-one
underapproximation for $M$, i.e., $V_{:1}W_{1:} \approx M (=R_1)$ and $V_{:1}W_{1:} \leq M$, we have that $R_2=M-V_{:1}W_{1:}$ is nonnegative. $R_2$ can then be underapproximated $V_{:2}W_{2:} \leq R_2$,
leading to $R_3 = R_2 - V_{:2}W_{2:}$, and so on. After $r$ steps, we get an underapproximation of rank $r$
\begin{eqnarray*}
M & \geq & V_{:1} W_{1:} + V_{:2} W_{2:} + \dots + V_{:r} W_{r:} \\
  &   =  & [V_{:1} \; V_{:2} \, \dots \, V_{:r}] [W_{1:}; \; W_{2:}; \, \dots \, ; W_{r:}] \\
  &   =  & V W.
\end{eqnarray*}

Besides enabling this recursive procedure, we notice that NMU leads to a more \emph{localized} part-based decomposition, in the sense that different basis elements tend to describe disjoint parts of the input data (i.e., involving different nonzero entries). This is a consequence of the underapproximation constraints which impose the extracted parts (the basis elements $V_{:k}$) to really be common features of the columns of $M$ since
\[
M_{:j} \; \gtrapprox \; \sum_{k} V_{:k} W_{kj}, \quad \forall j\;.
\]
Basis elements can only be combined to approximate a column of $M$ if each of them represents a part of this column, i.e., none of the parts selected with a positive weight can involve a nonzero entry corresponding to a zero entry in the input column $M_{:j}$. The following example demonstrates this behavior.

\begin{example}[Swimmer Database] \label{swimEx}
The swimmer image dataset consists of 256 binary images of a body with 4 limbs which can be each in 4 different positions. NMF is expected to find a part-based decomposition of these images, i.e., isolate different constitutive parts of the images (the body and the limbs) in each of its basis elements.

Figure~\ref{swim} displays a sample of such images along with the basis elements obtained with NMF and NMU. While NMF elements are rather sparse, they are mixtures of several limbs. By contrast, NMU returns a even sparser solution and is able to extract a single body part for each of its elements.
\begin{figure}[ht!]
\begin{center}
\includegraphics[width=\textwidth/2]{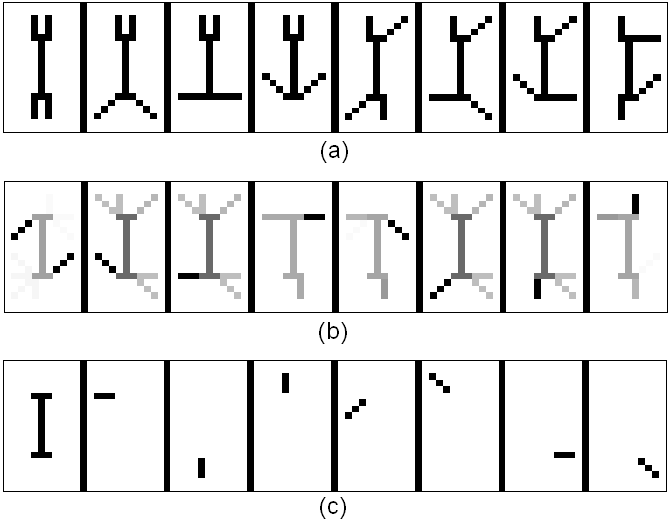}
\caption{Basis elements generated for the swimmer image dataset with $r=8$: (a) Sample images from the dataset, (b) NMF and (c) NMU; see Section~\ref{appl} for the algorithms used to compute the factorizations.}
\label{swim}
\end{center}
\end{figure}
\end{example}

\subsection{Sparsity}
The fact that NMU decompositions naturally generate sparser solutions than NMF can be explained as follows:
since the zero entries of $M$ can only be underapproximated by zeros, we have
\begin{equation*}  
M_{ij} = 0
\; \Rightarrow \; (VW)_{ij} = 0
\; \Rightarrow \;  V_{ik} = 0  \textrm{ or }  W_{kj} = 0, \, \forall k
\end{equation*}
which shows that when the input matrix is sparse, many components of the NMU factors will have to be equal to zero. This observation can be made more formal: defining the sparsity $s(M)$ of a $m$ by $n$ matrix $M$ as the proportion of its zero entries, i.e.,
\[
s(M) = \frac{\# \text{zeros}(M)}{mn} \; \in \; [0,1],
\]
we have the following theorem relating sparsity of $M$ and its NMU factors.
\begin{theorem} \label{sumsvsw}
For any nonnegative rank-one underapproximation $(v,w) \in \mathbb{R}^{m}_+ \times \mathbb{R}^{n}_+$ of $M \in \mathbb{R}^{m \times n}_+$ we have
\[
s(v) + s(w) \geq s(M).
\]
\end{theorem}
\begin{proof}
For a rank-one matrix $v w$, the number of nonzeros is exactly equal to the product of the number of nonzeros in vectors $v$ and $w$. Therefore we have that $(1-s(v w)) = (1-s(v)) (1-s(w))$ which implies $s(v w) = s(v) + s(w) - s(v)s(w) \le s(v) + s(w)$. Since underapproximation $v w$ satisfies $0 \leq vw \leq M$, it must have more zeros than $M$ and we have \[ s(M) \le s(sw) \le s(v) + s(w), \] proving our claim.
\end{proof}
Recall the recursive definition of the residuals $R_{k+1} = R_k - V_{:k}W_{k:}$ and $R_1 = M$. The following corollary relates their sparsity and the sparsity of the whole rank-$r$ approximation with that of the NMU factors.
\begin{corollary} \label{corosp}
For any nonnegative underapproximation $(V,W) \in \mathbb{R}^{m \times r}_+ \times \mathbb{R}^{r \times n}_+$ of $M \in \mathbb{R}^{m \times n}_+$ we have for each factor
\[
s(V_{:k}) + s(W_{k:}) \geq s(R_k) \geq s(M), \quad 1 \leq k \leq r,
\]
and $s(V) + s(W) \geq s(M)$.
\end{corollary}
\begin{proof}
We have $0 \leq V_{:k}W_{k:} \leq R_k \leq M$, which implies by the previous theorem the first set of inequalities. Observing that $s(V) = \frac{1}{r} \sum_k s(V_{:k})$ and $s(W) = \frac{1}{r} \sum_k s(W_{k:})$ is sufficient to prove the second inequality.
\end{proof}
Sparsity of the residuals $R_k$ is monotonically nondecreasing at each step, since $M = R_1 \geq R_2 \geq \dots \geq 0$. Moreover, the following theorem can guarantee an increase in sparsity at each step.
\begin{theorem} \label{zeroRes}
For any locally optimal nonnegative rank-one underapproximation $(v,w) \in \mathbb{R}^{m}_+ \times \mathbb{R}^{n}_+$ of $M \in \mathbb{R}^{m \times n}_+$, define sets $I$ and $J$ (supports of vectors $v$ and $w$) by \[
I = \{ i \;|\;  v_i > 0\}, \;
J = \{ j \;|\;  w_j > 0\},
\]
and define matrix $R(I,J)$ to be the submatrix of residual $R = M-vw$ whose row and column indices belong respectively to $I$ and $J$ (corresponding to the submatrix of $M$ that is not approximated by zeros). Then there is at least one zero in each row and each column of submatrix $R(I,J)$.
\begin{proof}
Simply observe that if $R(i,J) > 0$ (resp.\@ $R(I,j) > 0$) for some $i \in I$ (resp.\@ $j \in J$), $v_{i}$ (resp.\@ $w_{j}$) can be increased to obtain a strictly better solution, which contradicts the local optimality assumption.
\end{proof}
\end{theorem}

This ability of NMU to generate sparse part-based decomposition will be experimentally confirmed in Section~\ref{appl}.

\subsection{Related work}

The problem of rank-one underapproximation has been first introduced by Levin in \cite{Lev} in the case of \textsl{positive stochastic} matrices. He introduced a specific objective function different from the Frobenius norm and used a logarithmic change of variables in order to design an iterative method based on the corresponding optimality conditions.

In \cite{Gil}, the rank-one underapproximation problem is cast as a convex problem (hence efficiently solvable) using again different objective functions. Solutions are then used to initialize standard NMF algorithms in order to accelerate their convergence and, in general, find better final solutions as compared to those obtained with random initializations. Similar behavior was observed for other judicious initializations in~\cite{Alb, Bou, Cur}.

More recently, Dong et al.\@ \cite{DLC08} studied the same problem with the additional constraint that the rank of the residual must be strictly smaller than the rank of the factorized matrix. Using the Wedderburn rank reduction formula, they proposed a numerical procedure which is able to compute the maximum rank splitting of a nonnegative matrix. However, the underlying optimization problem is NP-hard \cite{Vav} and their algorithm is not guaranteed to find a solution in all cases.

Biggs et al.\@ \cite{Big} also introduced a recursive procedure to solve NMF problems: their idea is to locate and then approximate nearly rank-one submatrices of $M$. However, the problem of locating maximum rank-one submatrices is also shown to be NP-hard, and their algorithm is not globally optimal.

\subsection{Complexity}
\label{complex}

We now prove that \eqref{NMU} is NP-hard, even in the rank-one case (unlike \eqref{NMF}, which is polynomially solvable in the rank-one case). In order to do this, the rank-one version of the problem is proved to be equivalent to the biclique problem, which is NP-hard. The result is then generalized to \eqref{NMU} with arbitrary factorization rank $r$ using a simple construction.\\

A \emph{bipartite graph} $G_b$ is a graph whose vertices can be divided into two disjoint sets such that there is no edge between two vertices in the same set.
A \emph{biclique} $K_b$ is a complete bipartite graph, i.e., a bipartite graph where all the vertices from different sets are connected by an edge.
Finally, the so-called maximum edge biclique problem (the biclique problem for short) in a bipartite graph
\[
G_b = \Big(V = V_1 \cup V_2, E \subseteq (V_1 \times V_2) \Big),
\]
is the problem of finding a biclique $K_b = (V',E')$ in $G_b$ (i.e., $V'  = V'_1 \cup V'_2 \subseteq V$ and $E' = (V'_1 \times V'_2) \subseteq E$) with a maximum number of edges $|E'| = |V'_1|\cdot |V'_2|$.

Letting $M \in \{0,1\}^{m \times n}$ be the adjacency matrix of $G_b$ with $V_1 = \{s_1, \dots s_m\}$ and $V_2 = \{t_1, \dots t_n\}$, i.e.,
\[
M_{ij} = 1  \Leftrightarrow (s_i,t_j) \in E ,
\]
and introducing indicator binary variables $v_i$ (resp.\@ $w_j$) to denote whether $s_i$ (resp.\@ $t_j$) belongs to the biclique $K_b$, the Maximum edge Biclique Problem \eqref{MBP} in a bipartite graph can be formulated as follows
\begin{align}
\min_{v, w} \qquad & \sum_{i,j} (M_{ij}-v_iw_j)^2 \nonumber \\
&  v_i w_j \leq M_{ij}, \quad \forall i,j\;, \label{MBP}  \tag{MBP}  \\
&  v \in \{ 0,1 \}^{m}, \quad w \in \{ 0,1 \}^{n}. \nonumber
\end{align}
One can check  that this objective is equivalent to $\max_{v,w} \sum_{i,j} v_iw_j$. In fact, $M_{ij}-v_iw_j = (M_{ij}-v_iw_j)^2$ since $M$, $v$ and $w$ are binary and $M_{ij} \geq v_iw_j$.\\
The corresponding decision problem ``\textit{Given $K$, does $G_b$ contain a biclique with at least $K$ edges?}'' has been shown to be NP-complete~\cite{Peet}. Therefore, the corresponding optimization problem \eqref{MBP} is at least NP-hard.\\

For $r = 1$, \eqref{NMU} can be written as
\begin{align}
\min_{v \in \mathbb{R}^{m}, w \in \mathbb{R}^{n}} & \sum_{i,j} (M_{ij}-v_iw_j)^2 \nonumber \\
&  v_{i} w_{j} \leq M_{ij}, \;\, \forall i,j \;, \nonumber \label{NMU1}  \tag{NMU1}\\
& v \geq 0, \; w \geq 0, \nonumber
\end{align}
which is very close to \eqref{MBP}:
the difference is that vectors $v$ and $w$ are required to be \textit{binary} for \eqref{MBP} and \textit{nonnegative} for \eqref{NMU1}. The next lemma proves that the two problems are actually equivalent.

\begin{lemma}
\label{lemBin} For $M \in \{0,1\}^{m \times n}$, every optimal solution $(v,w)$ of \eqref{NMU1} is such that $v w$ is binary, i.e., $vw \in \{0,1\}^{m \times n}$, and can then be trivially transformed into a binary optimal solution $(v',w') \in \{0,1\}^m \times \{0,1\}^n$ of \eqref{MBP}.
\end{lemma}
\begin{proof} For $M=0$, this is trivial. Otherwise, suppose $(v,w)$ is an optimal solution of \eqref{NMU1}. Let define $(v',w')$ as
\begin{displaymath}
v_i' = \left\{ \begin{array}{ll}
1 & \textrm{if $v_i \neq 0$}\\
0 & \textrm{otherwise}
\end{array} \right.  \; \textrm{ and } \;
w_j' = \left\{ \begin{array}{ll}
1 & \textrm{if $w_j \neq 0$}\\
0 & \textrm{otherwise}
\end{array} \right. ,
\end{displaymath}
and analyze the different possibilities: as $v_i w_j \leq M_{ij}$, we have either
\begin{itemize}
\item $M_{ij} = 1$ and $0 < v_i w_j \leq 1 \Rightarrow v_i' w_j' = 1$;
\item $M_{ij} = 1$ and $v_i w_j = 0 \Rightarrow v_i' w_j' = 0$;
\item $M_{ij} = 0 \Rightarrow v_i w_j = 0 \Rightarrow v_i' w_j' = 0$.
\end{itemize}
Therefore, $v_i w_j \leq v_i' w_j' \leq M_{ij}$ which implies
\[
||M-v'w'||_F \leq ||M-vw||_F.
\]
By optimality of $(v,w)$, we must have $vw = v'w' \in \{0,1\}^{m \times n}$. Therefore, $(v',w') = (v/\max(v),w/\max(w))$ is an optimal binary solution of \eqref{NMU1} which is then also an optimal solution of \eqref{MBP} (note that we must have $\max(v) = \max(w)^{-1}$).
\end{proof}

\begin{corollary} \label{NMU1np}
\eqref{NMU1} is NP-hard.
\label{nmu1NP}
\end{corollary}

We now generalize Corollary~\ref{nmu1NP} to the more general case of \eqref{NMU} with $r>1$.
\begin{theorem}
\eqref{NMU} is NP-hard.
\label{nmuNP}
\end{theorem}
\begin{proof} Let $M \in \{0,1\}^{m \times n}$ be the adjacency matrix of a bipartite graph $G_b$. 
We define the matrix $A$ as
\begin{displaymath}
A = \textrm{diag}(M,r) = \left( \begin{array}{cccc}   M & 0 & \dots & 0 \\
0 & M & & 0 \\
\vdots & & \ddots & \vdots \\
0 & \dots & & M
 \end{array} \right) ,
\end{displaymath}
which is the adjacency matrix of another bipartite graph $G_b^r$ which is the graph $G_b$ repeated $r$ times.
Let $(V,W)$ be an optimal solution of \eqref{NMU}. Since $VW = \sum_{k=1}^r V_{:k} W_{k:}$, we have $VW \leq A \Rightarrow V_{:k}W_{k:} \leq A$. Therefore $(V_{:k},W_{k:})$ is a feasible solution of \eqref{NMU1} for the matrix $A$, i.e., for the graph $G_b^r$. Hence, each $(V_{:k},W_{k:})$ corresponds to a biclique $B_G^k =
(V_1^k \cup V_2^k, E^k)$ of $G_b^r$ with
\begin{displaymath}
V_{ik} \neq 0 \Leftrightarrow s_i \in V_1^k \quad \textrm{and} \quad W_{kj} \neq 0 \Leftrightarrow  t_j \in V_2^k.
\end{displaymath}
By optimality of $(V,W)$ and since there are at least $r$ independent maximum biclique in $G_b^r$, each $(V_{:k},W_{k:})$ must coincide with a maximum biclique of $G_b^r$ which corresponds to a maximum biclique of $G_b$. This is due to the fact that, because $G_b^r$ is the graph $G_b$ repeated $r$ times, a biclique clearly cannot span two disjoint subgraphs of $G_b^r$.
Therefore, \eqref{NMU} is NP-hard since any instance of \eqref{MBP} can be polynomially reduced to an instance of \eqref{NMU}.
\end{proof}

\section{An algorithm for NMU based on Lagrangian relaxation} \label{Lag}

Since \eqref{NMU}, like \eqref{NMF}, is a NP-hard problem, we can not expect to solve it up to guaranteed global optimality in a reasonable (e.g., polynomial) computational time (unless $P = NP$). In this section, we propose a nonlinear optimization scheme based on Lagrangian relaxation in order to compute approximate solutions of \eqref{NMU}.\\

\noindent Drop the $m \times n$ underapproximation constraints $VW \leq M$ of \eqref{NMU} and add them into the objective function with the corresponding Lagrange multipliers (dual variables, forming a matrix) $\Lambda \in \mathbb{R}^{m \times n}_+$, to obtain the Lagrangian function $L(V,W,\Lambda)$
\[
L(V,W,\Lambda) = \frac{1}{2}||M-VW||_F^2 + \sum_{i=1}^m \sum_{j=1}^n \Lambda_{ij} (VW-M)_{ij},
\]
where a factor of $\frac{1}{2}$ was introduced to make the presentation nicer. The Lagrangian relaxation subproblem \ref{lagrel} consists in minimizing $L$ for a fixed value of the $\Lambda$ multipliers, leading to the corresponding Lagrangian dual function $f(\Lambda)$
\begin{equation}
f(\Lambda) = \min_{V, W \geq 0} L(V,W,\Lambda) \label{lagrel} \tag{LR${}_\Lambda$}
\end{equation}
where $f(\Lambda)$ is well-defined because the minimum of $L(V,W,\Lambda)$ is always attained, due to the fact that $f$ is bounded below and the search space can be restricted to a compact set. Indeed, considering each rank-one factor individually $(V_{:k},W_{k:})$ and imposing w.l.o.g.\@ $||V_{:k}||_F^2 = ||W_{k:}||_F^2 = ||V_{:k}||_F ||W_{k:}||_F = ||V_{:k}W_{k:}||_F$, we have
\begin{eqnarray*}
||V_{:k}||_F^2 = ||W_{k:}||_F^2 
& \leq &  ||VW||_F \\
& \leq & ||M-\Lambda||_F + ||M-\Lambda-VW||_F,  \\
& \leq & 2||M-\Lambda||_F  \; \forall k
\end{eqnarray*}
where we have used the trivial solution $(V,W)=(0,0)$ to bound $||M-\Lambda-VW||_F$ (cf.\@ derivations of Section~\ref{VWup}).

Standard application of Lagrangian duality tells us that
\[
\text{(NMU)} \equiv \min_{V, W \geq 0} \quad \sup_{\Lambda \geq 0}  \; L(V,W,\Lambda) \quad \ge  \quad
\sup_{\Lambda \geq 0} \quad \min_{V, W \geq 0} \; L(V,W,\Lambda) \ = \ \sup_{\Lambda \geq 0} f(\Lambda),
\]
where the problem on the left of the inequality is equivalent to our original NMU formulation and the problem on the right is its Lagrangian dual, whose solution will provide a (hopefully tight) lower bound on the optimal \eqref{NMU}. This new problem is a nondifferentiable optimization problem with the nice property that its objective $f(\Lambda) = \min_{V, W \geq 0} L(V,W,\Lambda)$ is concave and its maximization (over a convex set) is then a convex problem (see \cite{AW07} and references therein).

We describe in the next section a general solution technique, which consists in repeatedly applying the following two steps: 
\begin{description}
\item[1.] Given multipliers $\Lambda$, compute $(V,W)$ to (approximately) minimize $L(V,W,\Lambda)$, i.e., solve \eqref{lagrel}; this is discussed in Section~\ref{VWup};
\item[2.] Given solution $(V,W)$, update multipliers $\Lambda$; this is described in  Section~\ref{Lamup}.
\end{description}

\subsection{Solving the Lagrangian relaxation problem}
\label{VWup}

The following derivations
\begin{align*}
L(V,W,\Lambda)
& =   \sum_{i,j} \frac{1}{2} (M-VW)_{ij}^2 + \sum_{i,j} \Lambda_{ij} (VW-M)_{ij}\\
& = \frac{1}{2} \sum_{i,j} M_{ij}^2 - \sum_{i,j}  M_{ij}(VW)_{ij} + \frac{1}{2} \sum_{i,j} (VW)_{ij}^2 \\
& \quad \;
+ \sum_{i,j} \Lambda_{ij} (VW)_{ij}
- \sum_{i,j} \Lambda_{ij} M_{ij}  \\
& =   \frac{1}{2} ||(M - \Lambda) - VW||_F^2 - \frac{1}{2} ||\Lambda||_F^2,
\end{align*}
show that minimizing $L(V,W,\Lambda)$ for a fixed $\Lambda$ is equivalent to minimizing
$||(M-\Lambda)-VW||_F^2$. Matrix $N=M-\Lambda$ is not necessarily nonnegative, therefore finding $V \geq 0$ and $W \geq 0$ such that \@ $N \approx VW$ is a more general problem than NMF. It is actually studied in detail  in~\cite{GG08} (see also \cite{Ding2}) where it is called Nonnegative Factorization~(NF), is formulated as
\begin{equation}
\min_{V \in \mathbb{R}^{m \times r}, W \in \mathbb{R}^{r \times n}}  ||N-VW||_F^2 \text{ such that }
 V\geq 0 \text{ and } W \geq 0, \label{NF} \tag{NF}
\end{equation}
with $N \in \mathbb{R}^{m \times n}$ and $1 \leq r < \min(m,n)$ and is shown to be NP-hard for any factorization rank (including $r=1$).

Some standard algorithms for NMF can easily adapted to handle an input matrix that is not nonnegative, i.e., solve a NF problem. For this work, we decided to use a recent technique called Hierarchical Alternating Least Squares (HALS), proposed in \cite{Cic}, which alternatively updates each column of $V$ and each row of $W$ with the following optimal closed-form solutions:
\begin{eqnarray}
V_{:k}^* & = & \textrm{argmin}_{V_{:k} \geq 0} \;\; ||(M-\Lambda)-VW||_F^2 \nonumber \\
& = & \max\Big({0}, \frac{A_{:k} - \sum_{l = 1, l \neq k}^{r} V_{:l} B_{lk}}{B_{kk}}\Big), \label{halsV}
\end{eqnarray}
with $A = (M-\Lambda)W^T$ and $B = WW^T$, and
\begin{eqnarray}
 W_{k:}^* & = & \textrm{argmin}_{W_{k:} \geq 0} \;\; ||(M-\Lambda)-VW||_F^2 \nonumber \\
 & = & \max\Big(\mathbf{0}, \frac{C_{k:} - \sum_{l = 1, l \neq k}^{r} D_{kl} W_{l:}}{D_{kk}}\Big), \label{halsW}
\end{eqnarray}
with $C = V^T(M-\Lambda)$ and $D = V^TV$. This can be viewed as a simple method of (block-)coordinate descent (also called alternating variables), which has been shown to perform strikingly well in practice, and much better than the popular multiplicative updates of Lee and Seung (see \cite{diep, Cic4, GG08}). Under some mild assumptions, every limit point of the above alternating scheme is a stationary point \cite{diep, GG08}.

The main computational cost of one HALS iteration is the evaluation of $A$ and $C$: they each require $2mnr$ (floating point) operations. One can check that the resulting total number of operations is $4mnr + O((m+n)r^2)$.

\begin{remark}
HALS is sensitive to the scaling of the initial matrices. For example, if the initial matrices $V$ and $W$ are chosen such that $VW \ggg M$, optimal columns of $V$ and optimal rows of $W$ computed by formulas \eqref{halsV} and \eqref{halsW} at the first step will most likely be equal to zero. This will lead to rank deficient approximations ($V_{:k}W_{k:} = 0$ for some $k$) and  numerical problems (for $V_{:k}=0$, update of $W_{k:}$ is not well defined and vice versa).
If the initial matrices $(V,W)$ are scaled~\cite{GG08}, i.e., by ensuring that
\begin{equation} \label{scaling}
1\, =\, \argmin_{\alpha} ||M-\alpha VW||_F = \frac{ \left\langle M,VW \right\rangle }{\left\langle VW,VW \right\rangle}
\end{equation}
where $\left\langle A,B\right\rangle = \sum_{i,j} A_{ij}B_{ij} = \mathrm{trace}(AB^T)$, this behavior is in general avoided. All initial matrices used in the following have been scaled.
\end{remark}

\subsection{Update of the multipliers $\Lambda$} \label{Lamup}

The second step of our algorithm consists in updating $\Lambda$ in order to find better (i.e., higher) solutions to the Lagrangian dual problem. Using the knowledge that any optimal solution $(\Lambda^*,V^*,W^*)$ of the Lagrangian dual must satisfy the following complementarity slackness conditions
\[
	\Lambda_{ij}^* (M-V^*W^*)_{ij} = 0 \ \forall i,j \text{, as well as feasibility conditions } \Lambda_{ij}^* \ge 0 \text{ and } (M-V^*W^*)_{ij} \ge 0 \text{, }
\]
we see that the update rule for the multipliers $\Lambda$ should satisfy the following: 
\begin{itemize}
\item if $(M-VW)_{ij} > 0$, $\Lambda_{ij}$ should be decreased and eventually reach zero if $(M-V^*W^*)_{ij} > 0$,
\item if $(M-VW)_{ij} < 0$, $\Lambda_{ij}$ should be increased to give more importance to $(M-VW)_{ij}$ in the cost
    function, hopefully in order to get a feasible solution such that $(M-V^*W^*)_{ij} \geq {0}$.
\end{itemize}
In the sequel, we use the following rule to update $\Lambda$, which satisfies the above requirements:
\begin{displaymath}
\Lambda \leftarrow  \max(0, \Lambda - \mu_k (M-VW)),\quad \mu_k \rightarrow 0,
\end{displaymath}
where $\mu_k$ is a predefined sequence of step lengths decreasing to zero; $\Lambda$ can be initialized to zero. This update is inspired from the concept of subgradient methods \cite{Shor85}; in fact, one can easily check that the quantity $(VW-M)$ is a subgradient of
\[
f(\Lambda) = \min_{V,W \geq 0} L(V,W,\Lambda)
\]
with respect to $\Lambda$, i.e., if $(\bar{V},\bar{W}) = \argmin_{V,W \geq 0} L(V,W,\bar{\Lambda})$, we have
\[
f(\Lambda)
\leq f(\bar{\Lambda})
     + \left\langle \bar{V}\bar{W}-M,{\Lambda}-\bar{\Lambda} \right\rangle
     , \quad \forall \Lambda.
\]
Two questions now arise
\begin{itemize}
\item Since an iterative algorithm is used to solve (approximately) the Lagrangian relaxation problem (cf.\@ section~\ref{VWup}), after how many of these HALS iterations do we stop and proceed to update the multipliers $\Lambda$?
\item How do we choose the sequence of step lengths $\mu_k$?
\end{itemize}
Subgradient methods usually assume that the Lagrangian relaxation problem~\eqref{lagrel} can be solved exactly and can guarantee their convergence to an optimal solution provided an appropriate sequence of step sizes is selected (see, e.g., \cite{AW07}), for example $\{\mu_k\}$ satisfying the conditions 
\[
0 \le \mu_k \rightarrow 0 \quad \textrm{ such that }  \quad  \sum_{k=0}^{\infty}  \mu_k^2 < +\infty \; \textrm{ while } \; \sum_{k=0}^{\infty}  \mu_k = +\infty.
\]
In the sequel, we choose to use $\mu_k = \frac{1}{k}$, which is such a suitable sequence. However, in our case, we cannot expect to solve \eqref{lagrel} in a reasonable amount of time since the problem is NP-hard. It would even probably be too expensive to wait for the stabilization of $(V,W)$ (e.g., getting close to a stationary but not necessarily optimal point). We therefore suggest to update $(V,W)$ only a constant number of times $T$ between each update of $\Lambda$, which leads to Algorithm L-NMU. Note that because we do not solve \eqref{lagrel} exactly, Algorithm L-NMU is not guaranteed to converge to an optimal solution of the Lagrangian dual but, as we will see, it produces satisfactory solutions in practice.
\algsetup{indent=2em}
\begin{algorithm}[ht!]
\caption{Lagrangian NMU (L-NMU)}  \label{LNMU}
\begin{algorithmic}[1]
\REQUIRE $M \in \mathbb{R}^{m \times n}_+$, $r > 0$, $V \in \mathbb{R}^{m \times r}_+$, $W \in \mathbb{R}^{r \times n}_+$, maxiter, $T$.
\ENSURE $(V,W)$ s.t. $VW \lesssim M$.\\
\medskip
\STATE $\Lambda = 0$; \\
\FOR {$k \, = \, 1$ : maxiter}
\STATE Update $(V,W)$ using $T$ iterations of HALS \eqref{halsV}-\eqref{halsW};
\STATE Update $\Lambda \leftarrow \max(0, \Lambda - \frac{1}{k} \, (M-VW))$;
\ENDFOR
\end{algorithmic}
\end{algorithm}

The additional computational cost of one iteration of algorithm L-NMU when compared with one iteration of HALS for NMF consists in the computation of $M-\Lambda$  (needed in step~3) and the update of $\Lambda$ (at step~4), which require $2mnr + O(mn)$ operations (and, in the special case $r=1$, $5mn$ operations).

\begin{remark}
Because convergence is not theoretically guaranteed, Algorithm L-NMU may end up with solutions that do not completely satisfy the underapproximation constraint. Although our numerical experiments show that this has no detrimental influence on the quality of the obtained sparse part-based representations (see Section~\ref{appl}), we give here a simple technique to transform such a solution into a feasible solution. Indeed, it is enough to consider the following QP problem  (convex quadratic objective function, linear inequality constraints) which only involves the $V$ factor
\begin{equation}
V^* = \argmin_{V \geq 0, VW \leq M} ||M-VW||_F^2. \label{NMUupdate}
\end{equation}
Because of its convexity, this problem can be solved up to global optimality in a very efficient manner, and replacing the original $V$ factor by the optimal solution $V^*$ leads to a feasible solution $(V^*,W)$ to \eqref{NMU}.
\end{remark}
\begin{remark}
Because update rule~\eqref{NMUupdate} is exact and computable in practice, it would be natural to consider a simpler algorithm based on its alternative application to the $V$ and $W$ factors,  without using the Lagrangian relaxation technique, hoping to converge to a solution of \eqref{NMU}. Unfortunately, we observed that this is quite inefficient in practice. In fact,
\begin{itemize}
\item it is relatively computationally expensive to solve these linearly constrained quadratic programs (with $mn+mr$ and $mn+nr$ inequalities), at least compared to the HALS closed-form update rules \eqref{halsV}-\eqref{halsW};
\item since the underapproximation constraint is imposed at each step, this algorithm has much less freedom to converge to good solutions: iterates rapidly get stuck on the boundary of the feasible domain, typically with (too) many zeros and a lower rank. For example, assuming $M$ has one zero in each column, we have that for any positive matrix $V$ the corresponding optimal $W$ is equal to 0:
\[
\forall j, \exists i \textrm{ s.t. } M_{ij} = 0
\Rightarrow \forall j, \exists i \textrm{ s.t. }
\sum_k V_{ik}W_{kj} = 0 \Rightarrow W_{kj} = 0, \forall k,j.
\]
Therefore, such an algorithm can only work if we decide a priori which values in $V$ and $W$ should be equal to zero, i.e. if we find a good sparsity pattern for the solution, which is precisely where the difficulty of the problem lies. 
Note that the same behavior is observed if a HALS-type algorithm is used instead of~\eqref{NMUupdate} (i.e., updating columns of $V$ and rows of $W$ alternatively): after the update of one column of $V$, the residual will have one zero in each row (cf.\@ Theorem~\ref{zeroRes}) which will prevent the other columns of $V$ to be nonzero (except if the sparsity pattern is chosen a priori).
\end{itemize}

\end{remark}

\begin{remark}
The L-NMU algorithm described above is not particularly well-suited to deal with very sparse input matrices. In fact, one has to store a potentially dense $m \times n$ matrix with the Lagrangian variables $\Lambda$. Nevertheless, Berry et al.\@ \cite{BGG} have obtained encouraging results when applying NMU to sparse anomaly detection problems in text mining. Moreover, it is possible to take advantage of the input sparsity pattern and design a computationally cheaper method. First note that the Lagrangian variables associated with a zero of $M$ will be nondecreasing in the course of the algorithm, since
\[
0 \leq (V^{(k)}W^{(k)})_{ij} \textrm{ and } (M)_{ij} = 0
\quad \Rightarrow \quad
\Lambda^{(k)}_{ij} \leq \Lambda^{(k+1)}_{ij},
\]
where superscript ${}^{(k)}$ denotes the solution at step $k$. 
Therefore one can significantly reduce the computational cost by defining
\[
\Lambda_{ij}^{(k)} = -g(k) \quad \textrm{ for all }  i,j \textrm{ s.t. } M_{ij} = 0,
\]
where $g(k)$ is an arbitrary positive nondecreasing function, e.g., $g(k) = \rho^k$ with $\rho > 1$, as is implicitly done in \cite{GG08}.

\end{remark}


\section{Numerical tests on image datasets} \label{appl}

We have argued in Section~\ref{nmuS} that NMU is potentially able to extract a better part-based representation of the data and that its factors should be sparser than those of the standard NMF, at the detriment of the approximation error. In this section, we support these claims by reporting results of computational experiments involving two variants of Algorithm L-NMU on several image datasets.

A direct comparison between NMU and NMF is not very informative in itself: while the former will provide a sparser part-based representation, the latter will feature a lower approximation error. This does not really tell us whether the improvements in the part-based representation and sparsity are worth the increase in approximation error. For that reason, we chose to compare NMU with two other sparse nonnegative matrix factorizations techniques, described below, in order to better assess whether the increase in sparsity achieved by NMU is worth the loss in reconstruction accuracy.

\subsection{Sparse NMF} \label{twosnmf}

We selected and tested the following two sparse nonnegative matrix factorization techniques that are frequently used in the literature.

\begin{itemize}
\item[1.] Hoyer describes in \cite{Hoy} an algorithm relying on additional explicit sparsity constraints on the factors, enforced at each iteration by means of a projection. The approximation error is reduced via a combination of projected gradient and multiplicative updates. For our experiments, we use the \matlab \, code provided by the author\footnote{This code was downloaded from \url{http://www.cs.helsinki.fi/u/phoyer/software.html}.}.

It should be pointed out that Hoyer is using a different definition of sparsity: for any nonzero $n$ dimensional vector $x$, his measure of sparsity ${sh}(x)$ is defined as
\begin{equation} \label{HoyS}
{sh}(x) = \frac{\sqrt{n} - ||x||_1/||x||_2}{\sqrt{n} - 1}
\quad \in \quad [0,1].
\end{equation}
Hence, a vector with a single nonzero entry is perfectly sparse
\[
\textrm{sh}([0 \dots 0 \; k \;  0 \dots 0]) = 1, \quad \forall k \neq 0,
\]
while a vector with all entries equal to each other is completely dense
\[
\textrm{sh}([k \dots k]) = 0,\quad \forall k \neq 0.
\]
In our experiments, we report sparsity using both the standard $s(\cdot)$ indicator and Hoyer's $sh(\cdot)$ measure.

\item[2.] Instead of enforcing sparsity at every iteration, a sparsity-inducing penalty term can be introduced in the objective function~\cite{LH01}. In particular, it is well-known that adding $l_1$-norm penalty terms induce sparser solutions (see, e.g., \cite{KP08, Kim, Cic4}), and we therefore solve the following problem:
\begin{equation}
\min_{V, W \geq 0} ||M-VW||_F^2 + \mu_V ||V||_1 + \mu_W ||W||_1,
\tag{sNMF} \label{sNMF}
\end{equation}
where $||A||_1 = \sum_{ik} |A_{ik}|$ and $\mu_V$ and $\mu_W$ are two positive parameters controlling the sparsity of $V$ and $W$. In order to solve \eqref{sNMF}, we use the HALS algorithm which can easily be adapted to handle the additional $l_1$-norm penalty terms (see, e.g., \cite{diep, Cic4}). This algorithm will be referred to as sNMF.
\end{itemize}

Technical details for the first technique are more complicated, but it allows the sparsity of the factors to be chosen a priori. The second technique is conceptually simpler but requires the determination of appropriate penalizing parameters by other means.


\subsection{Tested algorithms}

Algorithm~L-NMU proposed in Section~\ref{Lag} can be used to compute underapproximations for any given factorization rank $r$. This opens the possibility of building a rank-$r$ underapproximation in several different ways: one simple option consists in applying algorithm L-NMU directly to the rank-$r$ problem -- we call this method \textit{global NMU} (G-NMU). Another option consists in applying the recursive technique outlined in the introduction, used to motivate the introduction of underapproximations. More specifically, this means running algorithm L-NMU successively $r$ times to compute $r$ rank-one approximations, subtracting each approximation from the input matrix before computing the next one -- we call this method \textit{recursive NMU} (R-NMU). Note that many other variants are possible (e.g., computing two rank-$\frac{r}{2}$ approximations, computing $\frac{r}{2}$ successive rank-two approximations, etc.) but we only tested the two above-mentioned variants, which represent two extreme cases (no recursion and maximum recursion).

In both cases, our implementation of algorithm L-NMU computes two HALS steps between each update of the multipliers $\Lambda$ (i.e., we fixed $T = 2$). Most of the computational work done in one iteration of L-NMU consists in computing $M-\Lambda$, performing the two HALS steps and updating $\Lambda$; more specifically, one can estimate the computational cost of one iteration of G-NMU to $10mnr + O((m+n)r^2)$ operations, while an R-NMU iteration takes $13mn + O((m+n)r)$ operations (repeated $r$ times in the recursive procedure).

For each dataset, we test five nonnegative factorization algorithms: NMF based on HALS updates (NMF), global NMU (G-NMU), recursive NMU (R-NMU), sparse NMF with $l_1$-penalty terms (sNMF) and the algorithm of Hoyer. We also report the results of a standard Principal Component Analysis (PCA) to serve as a reference (recall that the approximation error of this unconstrained low-rank approximation, computed here with a singular value decomposition, is globally minimal for the given rank, but that its factors are neither nonnegative, nor sparse).

\subsection{Iteration limits and CPU time}

Each of the five iterative algorithms described above requires a limit on the number of its iterations; these limits were chosen in order to roughly allocate the same CPU time to each algorithm. More specifically, the standard NMF was given a 600-iterations limit, which corresponds to the computation of 600 HALS updates. The sparse sNMF, based on a slightly modified HALS update, was also allowed 600 iterations. Because a HALS update involves $4mnr + O((m+n)r^2)$ operations, we can deduce the following iteration budgets for G-NMU and R-NMU from the leading terms in their corresponding operation counts: G-NMU is allowed $600 \times \frac{4}{10} = 240$ L-NMU iterations while R-NMU can take $600 \times \frac{4}{13} \approx 180$ iterations.

An exception to the equal CPU time rule was made for the algorithm of Hoyer. Results obtained after an amount of CPU time similar to that of the other algorithms were too poor to be compared in a meaningful way. Indeed, because this method is based on a projected gradient method and multiplicative updates (both $O(mnr)$ operations per iteration), which are known to converge at a typically much slower rate, a relatively high limit of 1000 iterations had to be fixed, although the resulting CPU time is then much larger than for the other methods (for example, on the CBCL dataset, 600 iterations of HALS took $\sim 80s.$ while 1000 iterations of the algorithm of Hoyer needed $\sim 260s.$).

\subsection{Testing methodology}

Recall we decided to test algorithms sNMF and Hoyer to assess the quality of the sparsity-accuracy compromise proposed by our NMU approaches. To achieve this, we decided to pit each NMU variant against a solution of sNMF/Hoyer featuring the same level sparsity, and compare the resulting approximation errors. We therefore report results for $8$ algorithms on each dataset: PCA, NMF, G-NMU, sNMF with the same sparsity as G-NMU, which we denote by sNMF\{G-NMU\}, Hoyer\{G-NMU\}, R-NMU, sNMF\{R-NMU\} and Hoyer\{R-NMU\}.

In order to enforce a sparsity similar to the NMU solution in Hoyer's code, we compute the $sh$ measure of the NMU factors and input it as a parameter of the method (see description in subsection~\ref{twosnmf});  note however that we could only enforce this for the sparsest of the two NMU factors\footnote{Ideally, we would have imposed sparsity for both factors, but the implementation we used seemed to return poor results in that situation.}. In the case of sNMF, sparsity cannot be directly controlled, and penalty parameters are found using the following adaptive procedure, which proved to work well in practice: $\mu_V$ and $\mu_W$ are initialized to 0.1 and, after each iteration, $\mu_V$ (resp.\@ $\mu_W$) is increased by 5 percent if $s(V)$ (resp.\@ $s(W)$) is below the target sparsity, and is decreased by 5 percent otherwise.

All algorithms were run $10$ times with the same initial random matrices and only the best solution with respect to the Frobenius norm of the error is reported. When testing with gray-level images, the input matrices $M$ where normalized to have their entries varying between 0 and 1, with 0 representing white and 1 representing black (when trying to decompose $M$ as a sum of parts, this make more sense than the opposite convention, since the dark regions are the constitutive parts of the objects in the image datasets we analyze). Finally, before computing reported sparsity measures of the factors, any sufficiently small\footnote{We declare an entry of a factor to be sufficiently small if it is less than $0.1\%$ of the largest entry in its column.} entry is rounded to zero (indeed, because algorithms are stopped by the iteration limit before convergence, true zeros are typically not all reached). All tests were run within the \matlab \, 7.1 (R14) version, on a 3GHz Intel$^{\textrm{\textregistered}}$ Core{\texttrademark}2 Dual CPU PC.

\subsection{CBCL Face Database}

The CBCL face image dataset was used for the illustrative example of Figure~\ref{facto} and is made of 2429 gray-level images of faces represented with $19 \times 19$ pixels. We look for an approximation of rank $r=49$.

\begin{figure}[ht!]
\begin{center}
\includegraphics[width=\textwidth*3/4]{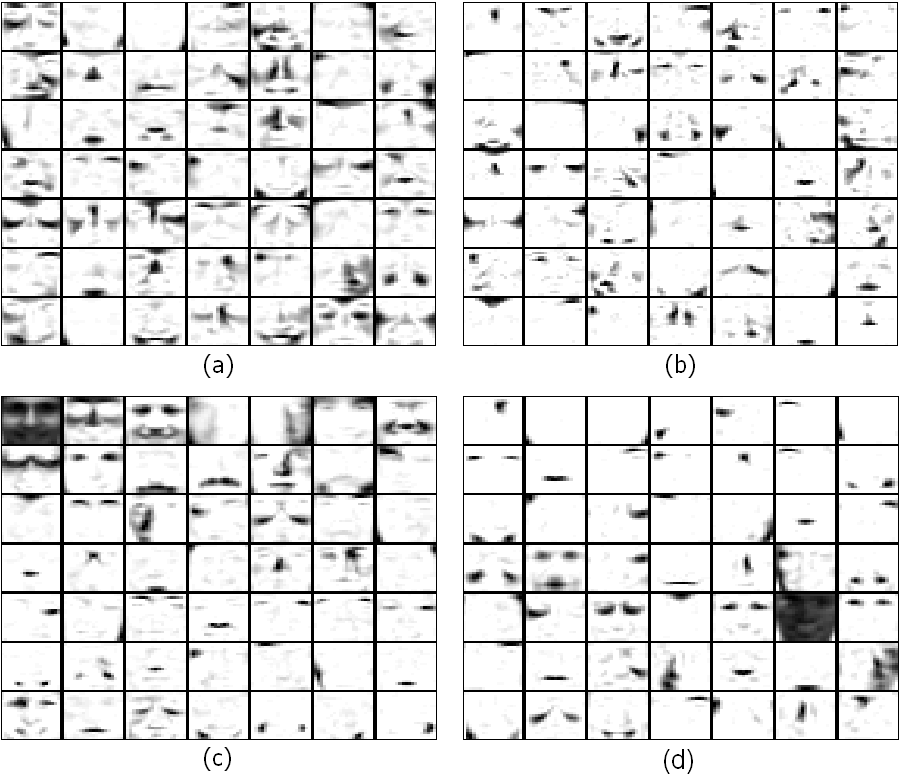}
\caption{Basis elements ($V_{:k}$) generated for the CBCL image dataset: (a)~NMF, (b)~G-NMU, (c)~R-NMU and (d)~sNMF with sparsity of G-NMU.}
\label{cbclB}
\end{center}
\end{figure}

Figure~\ref{cbclB} displays the basis elements for NMF, G-NMU, R-NMU and sNMF\{G-NMU\} (which was the best solution obtained in term of sparsity vs.\@ error among all four sNMF and Hoyer variants). Both G-NMU, R-NMU and sNMF achieve a better part-based representation than NMF, generating sparser solutions. An interesting feature of R-NMU is that it extracts parts successively in order of  ``importance'': the first basis element is a 'mean' face (which is dense) while the next ones describe different complementary parts (which become sparser as the recursion moves on, cf.\@ Corollary~\ref{corosp} and Theorem~\ref{zeroRes}).

A more quantitative assessment is provided in Table~\ref{cbclE}, reporting for the $8$ algorithms tested the relative error (in percent) of their solutions
\[
\textrm{ relative error} = \frac{||M-VW||_F}{||M||_F}
\]
in the second column (``Plain'') and the corresponding sparsity measures (in percent) of factors $V$ and $W$ in the last four columns.

\begin{table*}[ht!]
\begin{center}
\begin{tabular}{|c||c|c||c|c|c|c|c|}
\hline
\textit{}	 & Plain  & Improved &  $s(V)$ & $s(W)$ &  $sh(V)$ & $sh(W)$\\
\hline
PCA			  &	 7.43   &   	7.43 & 	 0	& 0 &  22  & 22 \\
\hline	
NMF			  &	 8.12   &   	8.11 & 	 56			 & 11 &  66  & 22 \\
\hline	
G-NMU &   12.45    &	   8.76 &  	74		 &  14 &  74  & 21 \\
sNMF\{G-NMU\} & 8.68   &   	8.44 & 	 74	  & 14 &  74  & 30 \\
Hoyer\{G-NMU\} & 9.33  &    8.78	 & 	69 	  & 6 &  73  & 16 \\
\hline
R-NMU	 &  16.42  &	10.89   &  53	 & 52 &   63 & 64 \\
sNMF\{R-NMU\}   &	 10.23    & 	9.49 &  	50		 & 50 &   56 & 57 \\
Hoyer\{R-NMU\} & 8.83  &   8.56	 & 	54 	  & 12 &  64  & 22 \\
\hline
\end{tabular}
\caption{Comparison of the relative approximation error and sparsity for the CBCL image dataset.}
\label{cbclE}
\end{center}
\end{table*}

As expected, PCA returns the smallest error, albeit with very dense factors. NMF already features much sparser factors (slightly half of the entries in $V$ are equal to zero), at the cost of a relatively modest increase in the approximation error ($7.43 \to 8.12$). G-NMU provides an even sparser solution (three quarters of zero entries), increasing again the approximation error ($8.12 \to 12.45$). The factors recursively computed by R-NMU are in comparison not as sparse: as explained above, this is because R-NMU focuses on obtained representative parts, including relatively dense ones for the first few steps of the recursion. However, it features a much sparse weight vectors, giving again more credit to the hypothesis that better parts are extracted. The corresponding approximation error is higher than for other methods, because the intrinsically greedy approach taken by R-NMU is not as efficient as a method that optimizes all the factors simultaneously.

Is the increased sparsity provided by G-NMU worth the increase in approximation error~? Looking at the corresponding results for sNMU\{G-NMU\} and Hoyer\{G-NMU\}, i.e., for sparse NMF and Hoyer's algorithms with a similar target sparsity, it might seem at first that the answer is negative: the other methods return solutions with similar number of nonzeros (slightly higher for Hoyer) and a lower approximation error ($8.68$ and $9.33$ instead of $12.45$). Actually, this was expected: because it tries to return an underapproximation, i.e., factors such that $V W \lesssim M$, the entries in the error term $M - VM$ are mostly nonnegative, while the other techniques, with no underapproximation constraint, obtain a smaller norm of the error by choosing the entries of $M-VW$ to be roughly half negative, half positive. It is therefore not completely fair to compare directly the error of the NMU approach to the other techniques.

In order to compensate for this, a simple rescaling could be used, i.e., multiplying $VW$ by a scalar since $VW \lesssim M$ (cf.\@ Equation~\eqref{scaling}). However, we chose a different procedure that has the advantage of benefiting all algorithms, including those whose error was not suffering from the underapproximation constraint. Once a solution is computed by one of the eight algorithms, we fix the zero entries of $V$ and $W$ and optimize the approximation error, i.e., $\min_{V,W \geq 0} ||M-VW||_F^2$, on the remaining (nonzero) entries (again, HALS can easily be adapted to handle this situation). In essence, this allows us to compare the sparsity patterns of the different solutions. We perform 100 additional HALS steps on each solution, and report the new relative approximation error in the third column of Table~\ref{cbclE} (``Improved''). Note that sNMF and Hoyer's errors are also improved by this procedure; this can be explained by the fact that they were also not directly trying to minimize the approximation error (Hoyer had to take into account its sparsity constraint, and sNMF was influenced by the penalty terms added to the approximation error).

Looking now at the NMU solutions in a fairer comparison, we observe that their approximation error becomes very close to that of sNMF and Hoyer, in particular for G-NMU, and not very far from the denser NMF, so that we can conclude that the sparsity-approximation error compromise it offers is worthwhile.

\subsection{Swimmer Database}

\begin{figure}[ht!]
\begin{center}
\includegraphics[width=\textwidth/2]{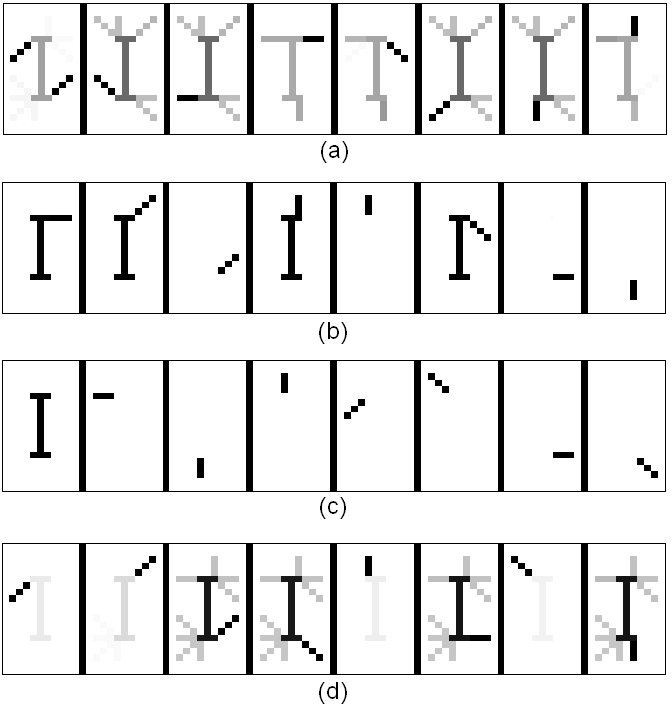}
\caption{Basis elements for the swimmer image dataset: (a)~NMF, (b)~G-NMU, (c)~R-NMU and (d)~sNMF with sparsity of R-NMU.}
\label{swimB}
\end{center}
\end{figure}

For the swimmer image dataset described in Example~\ref{swimEx} (256 images with $20 \times 11$ pixels), the $8$ basis elements obtained with the different algorithms are displayed on Figure~\ref{swimB} and the corresponding approximation errors and sparsity measures are reported in Table~\ref{swimE}.

\begin{table*}[ht!]
\begin{center}
\begin{tabular}{|c||c|c||c|c|c|c|}
\hline
\textit{Error}	 & Plain   & Improved &  $s(V)$ & $s(W)$ &  $sh(V)$ & $sh(W)$\\
\hline
PCA			  &	 37.98       & 	37.98 & 	 77*			 & 0 &  67  & 17 \\
\hline	
NMF			  &	 40.41       & 	40.41 & 	 84			 & 45 &  73  & 67 \\
\hline	
G-NMU &   47.70     &  46.85 &  	94		 &  75 &  85  & 78 \\
sNMF\{G-NMU\} & 50.52       & 	42.04 & 	 89	  & 66 &  84  & 73 \\
Hoyer\{G-NMU\} & 42.04      & 41.91	 & 	90 	  & 45 &  80  & 63 \\
\hline
R-NMU	 &  50.92     & 50.71   &  98	 & 66 &   93 & 65 \\
sNMF\{R-NMU\}   &	 41.66     & 	41.17 &  	85		 & 66 &   80 & 79 \\
Hoyer\{R-NMU\}**   &	 /     & 	/ &  	/		 &    / & / & /\\
\hline
\end{tabular}
\caption{Comparison of the relative approximation error and sparsity for the swimmer image dataset. *This value is very close to the percentage of zero rows in the matrix $M$ (corresponding to pixels that are equal to zero in all images): in general, PCA factors feature a zero component when all the entries of either one row or one column of the input matrix are equal to zero. **When imposing the sparsity level of R-NMU ($sh(V) = 0.93$), Hoyer's algorithm was not able to converge, probably because it is not well adapted to handle high sparsity constraints.
Note that sNMF is also sensitive to high sparsity requirements: high penalty terms sometimes lead to optimal zero factors ($V_{:k} = 0$ for some $k$), which had to be reinitialized.}
\label{swimE}
\end{center}
\end{table*}
As mentioned earlier, our two NMU algorithms are the only methods able to extract truly independent parts, while NMF and sNMF generate a combination of them. Note however that the solution generated by sNMF bears some similarity to the one of G-NMU.

\subsection{Hubble Space Telescope Spectral Images}

\begin{figure}[!ht]
\begin{center}
\includegraphics[width=\textwidth/2]{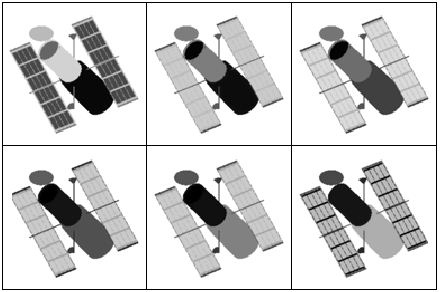}
\caption{Sample of Hubble space telescope spectral images.}
\label{hubble}
\end{center}
\end{figure}

\begin{figure}[ht!]
\begin{center}
\includegraphics[width=14cm]{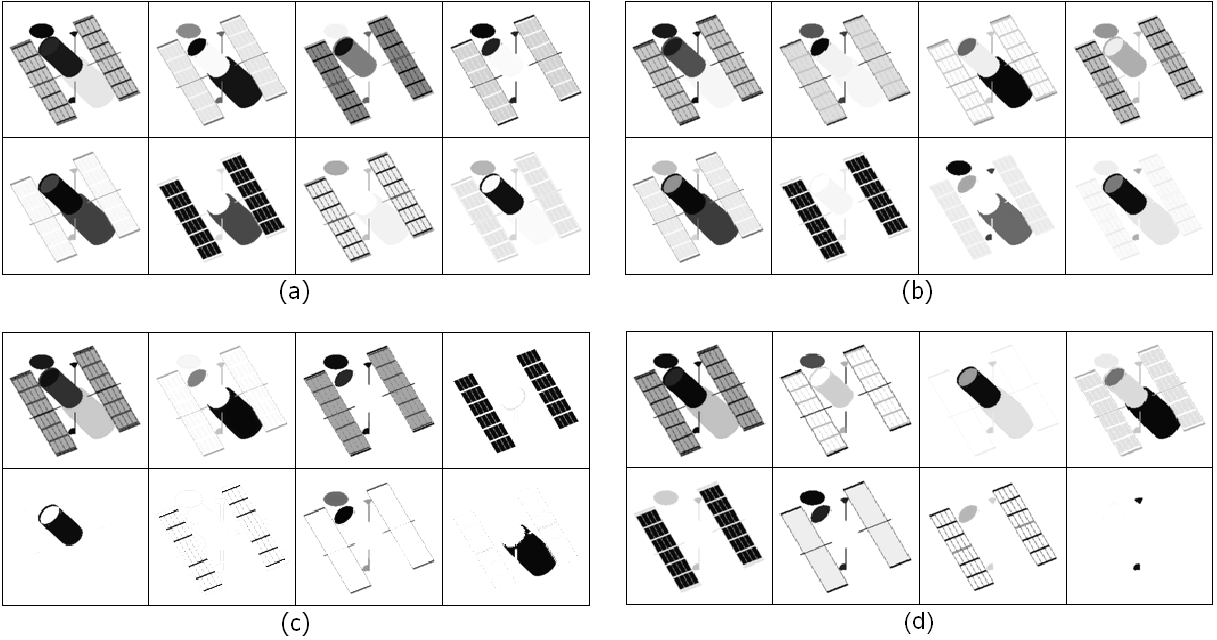}
\caption{Basis for the Hubble telescope: (a)~NMF, (b)~G-NMU, (c)~R-NMU and (d)~sNMF with sparsity of R-NMU.}
\label{hubbleB}
\end{center}
\end{figure}

The next image dataset consists of 100 spectral images ($128 \times 128$ pixels) of the Hubble telescope at different frequencies~\cite{hubble1, Plem}, see Figure~\ref{hubble}. With the choice $r=8$, NMF generates a nearly exact factorization (relative error $0.29\%$), because the spectral reflectance of the Hubble telescope results from the additive linear combination of the reflectance of eight constitutive materials. Figure~\ref{hubbleB} and Table~\ref{hubbleE} provide the visual and computational results for this dataset.

\begin{table*}[ht!]
\begin{center}
\begin{tabular}{|c||c|c||c|c|c|c|}
\hline
\textit{Error}	 & Plain  & Improved &  $s(V)$ & $s(W)$ &  $sh(V)$ & $sh(W)$ \\
\hline
PCA			  &	 0.01      & 	0.01 & 	 57			 & 0 &   62  & 25\\
\hline	
NMF			  &	 0.29     & 	0.29 & 	 64			 & 5 &   57  & 35\\
\hline	
G-NMU &   0.52      &  0.11 &  	64		 &  4 &   60   & 31\\
\hline
R-NMU	 &  3.74     &  1.37   &  79	 & 30 &   71  & 62 \\
sNMF\{R-NMU\}   &	 0.48     & 	0.37 &  	73		 & 28 &  66   & 64 \\
Hoyer\{R-NMU\} &   0.77      & 	0.68 & 	75 	  & 0 &  71  & 12 \\
\hline
\end{tabular}
\caption{Comparison of the relative approximation error and sparsity for the Hubble telescope image dataset.}
\label{hubbleE}
\end{center}
\end{table*}
Because NMF is already a nearly exact reconstruction  (Table~\ref{hubbleE}), the NMU constraints are somehow redundant: NMF and G-NMU are basically equivalent and return solutions with very similar sparsity measures (albeit with a slightly lower error for G-NMU). For that reason, sNMF\{G-NMU\} and Hoyer\{G-NMU\} return results nearly identical to NMF and are omitted from the table.

Recursive R-NMU extracts parts in order of importance: first, a global picture of the telescope and then its different constitutive parts. This allows it to generate the sparsest solution, with several basis elements representing well-delimited constitutive parts of the telescope not identified by the other methods.

%

\subsection{Kuls Illuminated Faces}

\begin{figure*}[ht!]
\begin{center}
\includegraphics[width=\textwidth*2/3]{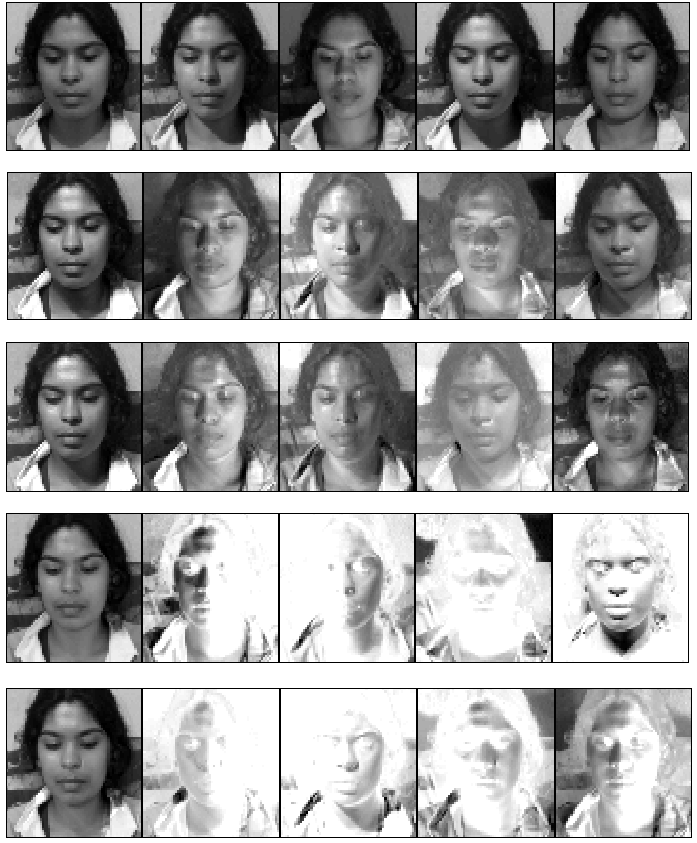}
\caption{Basis for the Kuls image dataset, from top to bottom: sample of images, NMF, G-NMU, R-NMU, sNMF with sparsity of R-NMU.}
\label{Basiskuls}
\end{center}
\end{figure*}

A static scene was illuminated from many directions with a moving light source to produce the Kuls image dataset\footnote{Available at \url{http://www.robots.ox.ac.uk/~amb/}.}. It consists of 20 images ($64 \times 64$ pixels) of a face. Because the images are very similar, most of the information (more than 70 percent) can be expressed with only one factor. The remaining information resides in the different orientations of the lighting.
Computational and visual results for a rank-5 factorization are given by Table~\ref{kulsE} and Figure~\ref{Basiskuls}.
We observe that NMF and G-NMU obtain similar results: even though they are both able to extract several faces with different lighting orientations, they do not extract a sparse and part-based representation.

R-NMU first extracts a face illuminated from all directions, and then complementary parts representing different orientations of the lighting (successively on the fourth row of Figure~\ref{Basiskuls}: global then light from the right, left, bottom and top). This nice recursive extraction of the information is a direct consequence of the underapproximation constraints. Although sNMF (with the same sparsity requirement as R-NMU) is also able to extract a part-based representation with a slightly better approximation error, only two components are well-identified (left and right lighting mixed with top and bottom lighting).
\begin{table*}[ht!]
\begin{center}
\begin{tabular}{|c|c|c|c|c|c|c|c|c|}
\hline
\emph{Error}	 & Plain  & Scaled & Improved &&  $s(V)$ & $s(W)$ &  $sh(V)$ & $sh(W)$\\
\hline
PCA			  &	 4.36   &   4.36   & 4.36 &&  	 	0		 &  0 &  23   & 15\\
\hline
NMF			  &	 4.38   &   4.38   & 4.38 &&  	 	1		 &  7 &  9   & 38\\
\hline	
G-NMU &    6.27   &	 5.77   &  4.49 &&   	3		 & 20 &   8  & 48\\
sNMF\{G-NMU\}  & 4.42  &   4.42   & 4.41	 && 	 	2  & 20 &  8   & 47\\
Hoyer\{G-NMU\}   & 4.60  &    4.60   & 4.71	 && 	 	2  & 25 &  8   & 53\\
\hline
R-NMU	 &  8.13  &	 7.84   &  5.73  &&  	29 &  31  &  38   & 67\\
sNMF\{R-NMU\}   &	 5.24  &  5.24  & 5.01	 &&  		29	 &  31  &  32    & 59\\
Hoyer\{R-NMU\}  & 6.82  &   6.82   & 6.54	 &&	 	0  & 71 &  6   & 92\\
\hline
\end{tabular}
\caption{Comparison of the relative approximation error and sparsity for the Kuls image dataset.}
\label{kulsE}
\end{center}
\end{table*}

\section{Conclusion}

In order to solve the NMF problem in a recursive way, we have introduced a new
problem, namely Nonnegative Matrix Underapproximation (NMU), which was shown
to be NP-hard using its equivalence with the maximum-edge biclique problem. The additional constraints of NMU are shown to induce sparser factors and to lead naturally to a better part-based representation of the data, while keeping a fairly good reconstruction. We proposed an algorithm based on Lagrangian relaxation to find approximate solutions to NMU.


We tested two factorization methods based on this algorithm, one with full recursion (R-NMU), the other without recursion (G-NMU), on several standard image datasets. After suitable post-processing, we observed that the factors computed by these methods indeed offer a good compromise between their achieved sparsity and the resulting approximation error, comparable or sometimes superior to that of two standard sparse matrix factorization techniques.

These two variants can be contrasted in the following way: where G-NMU mainly focuses on finding sparse
factors with small reconstruction error, in the same spirit as sNMF and Hoyer, R-NMU typically computes an even sparser factorization corresponding to a better part-based representation, albeit with a moderate increase in the
reconstruction error (due to the greedy approach). Moreover, this second variant is useful in situations where the factorization rank is not fixed a priori: the fact that it is recursive allows the user to stop the procedure as soon as the reconstruction error becomes satisfactory, without having to recompute a completely different solution from scratch every time a higher-rank factorization needs to be considered.

\section*{Acknowledgments}
The authors would like to thank the anonymous reviewers for their insightful comments which helped improve the paper.

\end{document}